\documentclass{article}
\usepackage{eurosym}
\usepackage{amssymb}
\usepackage{amsfonts}
\usepackage{amsmath}

\begin{document}

\begin{center}
{\Large Some properties of the norm in a quaternion division algebra }

\begin{equation*}
\end{equation*}%
Cristina FLAUT and Diana SAVIN 
\begin{equation*}
\end{equation*}
\end{center}

\textbf{Abstract. }{\small In this paper we provide some applications of the
norm form in \ some quaternion division algebras over rational field and we
give some properties of Fibonacci sequence and Fibonacci sequence in
connection with quaternion elements. We define a monoid structure over a
finite set on which we will prove that the defined Fibonacci sequence is
stationary, we provide some properties of the norm of a rational quaternion
algebra, in connection to the famous Lagrange's four-square theorem and its
generalizations given by Ramanujan. Moreover, we prove some results
regarding the arithmetic of integer quaternions defined on some division
quaternion algebras and we define and give properties of some special
quaternions by using Fibonacci sequences. }%
\begin{equation*}
\end{equation*}


\medskip

\begin{equation*}
\end{equation*}

\bigskip \textbf{1. Introduction}%
\begin{equation*}
\end{equation*}

Over the time, the problem of determining integer numbers $a,b,c,d$ such
that all positive integers can be written under the form $%
ax^{2}+by^{2}+cz^{2}+du^{2}$, with $x,y,z,u$$\in $$\mathbb{Z}$, has
fascinated many mathematicians. Jacobi proved that all positive integers $n$
can be written under the form $n=x^{2}+2y^{2}+3z^{2}+6u^{2}$, with $x,y,z,u$$%
\in $$\mathbb{Z}$. Liouville and Pepin proved that for 
\begin{equation*}
\left( a,b,c,d\right) \in \left\{ \left( 1,1,1,1\right) ,\left(
1,1,2,2\right) ,\left( 1,2,2,2\right) ,\left( 1,1,1,4\right) ,\left(
1,1,2,4\right) ,\left( 1,2,2,4\right) \right\}
\end{equation*}%
\begin{equation*}
\cup \left\{ \left( 1,2,4,4\right) ,\left( 1,1,2,8\right) ,\left(
1,2,4,8\right) ,\left( 1,1,3,3\right) ,\left( 1,2,5,10\right) \right\} ,
\end{equation*}%
all positive integers can be written under the form $%
n=ax^{2}+by^{2}+cz^{2}+du^{2}$, with $x,y,z,u$$\in $$\mathbb{Z}.$\newline

In the paper [Ra; 1917], S. Ramanujan has found all positive integers $%
a,b,c,d$ such that all positive integers $n$ can be expressed under the form 
\begin{equation*}
ax^{2}+by^{2}+cz^{2}+du^{2},
\end{equation*}
with $x,y,z,u$$\in $$\mathbb{Z}.$ He proved that there are only $55$ sets of
values of $a,b,c,d$ such that any positive integer $n$ can has the form $%
n=ax^{2}+by^{2}+cz^{2}+du^{2}$, with $x,y,z,u$$\in $$\mathbb{Z}$ ([Ra;
1917], p. 171). From these $55$ sets of values, only $7$ sets of values of $%
a,b,c,d~$ have the properties $a=1$ and $d=bc$, namely for 
\begin{equation}
\left( a,b,c,d\right) \in \left\{ 
\begin{array}{c}
\left( 1,1,1,1\right) ,\left( 1,1,2,2\right) ,\left( 1,1,3,3\right) ,\left(
1,2,2,4\right) ,\left( 1,2,3,6\right) , \\ 
\left( 1,2,4,8\right) ,\left( 1,2,5,10\right)%
\end{array}%
\right\} .  \tag{1.1.}
\end{equation}%
\qquad

In this paper we find an interesting connection between quaternion algebras
and these quaternary quadratic forms. Moreover, we define a monoid structure
over a finite set on which we will prove that the defined Fibonacci sequence
is stationary, we prove some results regarding the arithmetic of integer
quaternions defined on some division quaternion algebras and we define and
give properties of some special quaternions by using Fibonacci sequences.%
\begin{equation*}
\end{equation*}

\textbf{2.} \textbf{Fibonacci stationary elements}%
\begin{equation*}
\end{equation*}

\bigskip Let $\left( f_{n}\right) _{n\geq 0}$ be the Fibonacci sequence 
\begin{equation*}
f_{n}=f_{n-1}+f_{n-2},\;n\geq 2,f_{0}=0;f_{1}=1,
\end{equation*}%
and $\left( l_{n}\right) _{n\geq 0}$ be the Lucas sequence 
\begin{equation*}
l_{n}=l_{n-1}+l_{n-2},\;n\geq 2,l_{0}=2;l_{1}=1.
\end{equation*}%
We know that the Binet's formula for Fibonacci sequence are 
\begin{equation*}
f_{n}=\frac{\alpha ^{n}-\beta ^{n}}{\alpha -\beta }=\frac{\alpha ^{n}-\beta
^{n}}{\sqrt{5}},\ \ \left( \forall \right) n\in \mathbb{N}.
\end{equation*}%
and the Binet's formula for Lucas sequence are 
\begin{equation*}
l_{n}=\alpha ^{n}+\beta ^{n},\ \ \left( \forall \right) n\in \mathbb{N}.
\end{equation*}

\textbf{Proposition 2.1.} ([Fib]). \textit{Let} $(f_{n})_{n\geq 0}$ \textit{%
be the Fibonacci sequence} \textit{and let } $(l_{n})_{n\geq 0}$ \textit{be
the Lucas sequence.} \textit{The following properties hold:}\newline
i) 
\begin{equation*}
f_{n}^{2}+f_{n+1}^{2}=f_{2n+1},\forall ~n\in \mathbb{N};
\end{equation*}%
ii) 
\begin{equation*}
f_{n+3}^{2}=2f_{n+2}^{2}+2f_{n+1}^{2}-f_{n}^{2},\forall ~n\in \mathbb{N}%
-\{0\};
\end{equation*}%
iii) 
\begin{equation*}
f_{n-1}f_{m}+f_{n}f_{m+1}=f_{n+m},\forall ~n,m\in \mathbb{N}-\{0\};
\end{equation*}%
iv) 
\begin{equation*}
f_{n}f_{m}-f_{n-k}f_{m+k}=\left( -1\right) ^{n-k}f_{k}f_{m+k-n},\forall
~n,m,k\in \mathbb{N},k\leq n.
\end{equation*}

It is well known that the Fibonacci sequence is periodic when it is reduced
modulo $m$ and this period is called \textit{Pisano's} period, denoted by $%
\pi \left( m\right) $. Fibonacci numbers are defined by using the addition
operation over the real field. Will be interesting to find other finite
algebraic structures on which we will study the behavior of such a sequence.
Therefore, it is interesting to find other multiplications and sets such
that the behavior of such a sequence defined on to be predictable. In the
following, we will define a monoid structure over a finite set on which we
will prove that the defined Fibonacci sequence is stationary.

Let $\left( M_{1},\perp \right) $ and $\left( M_{2},\intercal \right) $ be
two monoids. On their cartesian product $M_{1}\times M_{2}$ we define the
following multiplication $"\ast "$:%
\begin{equation}
\left( a_{1},a_{2}\right) \ast \left( b_{1},b_{2}\right) =\left( a_{1}\perp
b_{1},a_{2}\intercal b_{2}\right) .  \tag{2.1.}
\end{equation}%
It results that $\left( M_{1}\times M_{2},\ast \right) $ is a monoid.

Now, we consider the totally ordered sets $\left( A_{i},\leq _{i}\right) $,
with $A_{i}=\{o_{i},x_{i}\},i\in \{1,2,...,k\}$. \ On $A_{i}$ we define the
multiplication $"+_{i}"$%
\begin{equation*}
x_{i}+_{i}x_{i}=x_{i},o_{i}+_{i}o_{i}=o_{i},x_{i}+_{i}o_{i}=o_{i}+_{i}x_{i}=x_{i}.
\end{equation*}%
It results that $\left( A_{i},+_{i}\right) $ is a commutative monoid. We
consider the Cartesian product $A=A_{1}\times A_{2}\times ...\times A_{k}.$
Using relation $\left( 2.1\right) $, we obtain that $\left( A,\ast \right) $
is a commutative monoid. The set $A=\{y_{0}\leq y_{1}\leq ...\leq
y_{2^{k}-1}\}$ is ordered with lexicographic order and it is a totally
ordered set. We denote $y_{0}=0=\left( o_{1},o_{2},...,o_{2^{k}}\right) $
and $y_{2^{k}-1}=1=\left( x_{1},x_{2},...,x_{2^{k}}\right) $.

On $A$, the above multiplication $"\ast "$ can be given by the following
relations:

\begin{equation*}
\left\{ 
\begin{array}{c}
y_{i}\ast y_{j}=1\text{, if }i+j>2^{k}-1\text{;} \\ 
y_{i}\ast y_{j}=y_{i+j}\text{, if }i+j\leq 2^{k}-1\text{;} \\ 
y\ast 0=y,\text{ for all }y\in A.%
\end{array}%
\right. .
\end{equation*}

\textbf{Example 2.2.} Let $k=3$. We have $A=\{\left(
y_{0},y_{1},...,y_{7}\right) \}$, where\newline
$y_{0}=\left( o_{1},o_{2},o_{3}\right) =0,y_{1}=\left(
o_{1},o_{2},x_{3}\right) ,y_{2}=\left( o_{1},x_{2},0_{3}\right)
,y_{3}=\left( o_{1},x_{2},x_{3}\right) ,$\newline
$y_{4}=\left( x_{1},o_{2},o_{3}\right) ,y_{5}=\left(
x_{1},o_{2},x_{3}\right) ,y_{6}=\left( x_{1},x_{2},o_{3}\right)
,y_{7}=\left( x_{1},x_{2},x_{3}\right) =1$. It results that

$y_{1}\ast y_{2}=\left( o_{1},x_{2},x_{3}\right) =y_{3},y_{2}\ast
y_{5}=\left( x_{1},x_{2},x_{3}\right) =y_{7}=1,$

$y_{3}\ast y_{6}=\left( x_{1},x_{2},x_{3}\right) =1,$ etc.\medskip

\textbf{Definition 2.3. }We consider $\left( M,\ast \right) $ a commutative
monoid. For $a,b\in M$, we define the following sequence 
\begin{equation*}
\{a,b\}=\{a,b,a\ast b,b\ast \left( a\ast b\right)
,...,.v_{n},v_{n+1},v_{n+2},.....\},
\end{equation*}%
where $v_{0}=a,v_{1}=b$ and $v_{n+2}=v_{n}\ast v_{n+1}$, for $n\in \mathbb{N}
$. This sequence are called the \textit{Fibonacci sequence attached to the
elements} $a,b$. If we find a number $t\in \mathbb{N}$ such that $%
v_{n}=v_{n+1}=v_{n+2}=....$, for all $n\geq t$, then the sequence $\{a,b\}$
is called $t$\textit{-stationary}.\medskip

\textbf{Proposition 2.4.} \textit{With the above notations, the Fibonacci
sequence defined on the set} $A$ \textit{is} $t$\textit{-stationary.\medskip 
}

\textbf{Proof.} If $a,b\in M,a\neq 0,b\neq 0,$ we consider the sequence%
\newline
\begin{equation*}
\{a,b\}=\{a,b,a\ast b,b\ast \left( a\ast b\right)
,...,v_{n},v_{n+1},v_{n+2},.....\},
\end{equation*}%
where $v_{0}=a,v_{1}=b$ and $v_{n+2}=v_{n}\ast v_{n+1}$, for $n\in \mathbb{N}
$. It results that $v_{2}=a\ast b>a$ and $v_{2}=a\ast b>b$, since $a\ast b=a 
$ or $a\ast b=b$ implies $b=0$ or $a=0$, false. If $v_{2}=a\ast b>ax$, we
have that\newline
$v_{3}=b\ast \left( a\ast b\right) =\left( a\ast b\right) \ast b>a\ast
b=v_{2},$\newline
$v_{4}=(a\ast b)\ast \left( b\ast \left( a\ast b\right) \right) >b\ast
\left( a\ast b\right) =v_{3}$, etc. We get that the obtained increased
sequence is stationary, since the set $A$ is finite. Therefore, there is
number $t\in \mathbb{N}$ such that $v_{t}=1=v_{t+1}=v_{t+2}=...$. It results
that

\begin{equation*}
\{a,b\}=\{a,b,a\ast b,b\ast \left( a\ast b\right) ,...,v_{t-1},1,1,.....\}.
\end{equation*}

Assuming now that $a,b\in A$ and $b=0$, we obtain the sequence $\{a,b\}$
with $v_{0}=a,v_{1}=0,v_{2}=a,v_{3}=a,v_{4}=$ $a\ast a>a$. Therefore, from
the above, we obtain a stationary sequence, that means we get a number $t\in 
\mathbb{N}$ such that $v_{t}=1=v_{t+1}=v_{t+2}=...$.$_{{}}\medskip $

\textbf{Example 2.5.} Considering the set $A$ from Example 2.2, and the
elements $y_{2}=\left( o_{1},x_{2},o_{3}\right) ,y_{4}=\left(
x_{1},o_{2},o_{3}\right) $, we have\newline
$\{y_{2},y_{4}\}=\{y_{2},y_{4},y_{2}\ast y_{4},y_{4}\ast \left( y_{2}\ast
y_{4}\right) ,....\}=$\newline
$=\{y_{2},y_{4},y_{6},1,1,....\}$, therefore the sequence is
stationary.\medskip 

\textbf{Remark 2.6.} Proposition 2.4 from above generalizes results obtained
in \ [Fl; 20], Proposition 2.6.

\begin{equation*}
\end{equation*}

\bigskip \textbf{3. Remarks regarding the norm of a quaternion algebra} 
\begin{equation*}
\end{equation*}

In this section, we will provide some properties of the norm of a rational
quaternion algebra, in connection to the famous Lagrange's four-square
theorem and its generalizations given by Ramanujan in \textbf{[}Ra; 1917%
\textbf{].}

First of all, we recall some results about diophantine equations, about
quaternion algebras and about Fibonacci sequence, results which will be used
in the following proofs.\medskip\ \medskip \newline

\textbf{Proposition 3.1.} ([Si; 70]). \textit{Let} $m$ \textit{be a fixed
positive integer. The diophantine equation} $x^{2}+my^{2}=z^{2}$ \textit{has
infinitely many solutions, namely:} 
\begin{equation*}
x=a^{2}-mb^{2},y=2ab,z=a^{2}+mb^{2},a,b\in \mathbb{Z}.
\end{equation*}%
$_{{}}\medskip $

\textbf{Proposition 3.2.} ([Gi, Sz; 06]). \textit{Let} $K$ \textit{be a
field. Then, the quaternion algebra} $\mathbb{H}_{K}\left( b,c\right) $ 
\textit{is a split algebra if and only if the conic} 
\begin{equation*}
C\left( a,b\right) :bx^{2}+cy^{2}=z^{2}
\end{equation*}%
\textit{has a rational point over} $K$$\backslash $ $\left\{ \left(
0,0,0\right) \right\} $, \textit{i.e. there are} $x_{0},y_{0},z_{0}\in K$$%
\backslash $$\left\{ \left( 0,0,0\right) \right\} $ \textit{such that} $%
bx_{0}^{2}+cy_{0}^{2}=z_{0}^{2}._{{}}\medskip $\newline

Let $K$ be an algebraic number field, let $\mathcal{O}_{K}$ be the ring of
integers of the field $K$ and let $H_{K}(b,c)$ be a quaternion algebra over
the field $K.$ We recall that \textit{the discriminant (or reduced
discriminant)} $D_{H_{K}(b,c)}$ of the quaternion algebra ${H}_{K}(b,c)$ is
the product of all prime ideals of $\mathcal{O}_{K}$ that ramify in ${H}%
_{K}(b,c).$ If $p$ is a prime and $K=\mathbb{Q},$ we recall the notations: $%
\left( \frac{\cdot }{p}\right) $ for the Legendre symbol in $\mathbb{Z}$ and 
$\left( b,c\right) _{p}$ for the Hilbert symbol in the p-adic field $\mathbb{%
Q}_{p}$ 
\begin{equation*}
\left( b,c\right) _{p}=\left\{ 
\begin{array}{c}
1,\text{if}~p\ \text{does\ not\ ramify\ in}\ {H}_{K}(b,c) \\ 
-1,\text{if}~p\ \text{ramifies\ in}\ {H}_{K}(b,c).%
\end{array}%
\right. .
\end{equation*}%
\newline
\qquad The quaternion algebra $H_{K}(b,c)$ splits if and only if the
discriminant $D_{H_{K}(b,c)}$ is equal to $\mathcal{O}_{K}.$ Since a
quaternion algebra is either split or a division algebra, it results that a
quaternion algebra ${H}_{K}(b,c)$ is a division algebra if and only if there
is a prime $p$ such that $p|D_{H_{K}(b,c)}$.\medskip \newline
\qquad \textbf{Proposition 3.3.} ([Ma; 07]; [Vi; 80]). \textit{The
discriminant of a quaternion algebra determines the algebra up to
isomorphism. That means, if} $K$ \textit{is a field, then two quaternion
algebras over the field} $K$ \textit{are isomorphic if and only if they have
the same reduced discriminant.}$_{{}}\medskip $

Let $b,c$ be positive integer numbers and let $\mathbb{H}_{\mathbb{Q}}\left(
-b,-c\right) $ be the generalized quaternion algebra over the field of
rational numbers. Let $x$$\in $$\mathbb{H}_{\mathbb{Q}}\left( -b,-c\right) ,$
$x=x_{1}\cdot 1+x_{2}e_{2}+x_{3}e_{3}+x_{4}e_{4},$ where $x_{i}\in \mathbb{Q}%
,i\in \{1,2,3,4\}$, and the elements of the basis $\{1,e_{2},e_{3},e_{4}\}$
satisfy the following rules, given in the below multiplication table:

\begin{center}
\begin{tabular}{ccccc}
$\cdot $ & $1$ & $e_{2}$ & $e_{3}$ & $e_{4}$ \\ \hline
$1$ & $1$ & $e_{2}$ & $e_{3}$ & $e_{4}$ \\ 
$e_{2}$ & $e_{2}$ & $-b$ & $e_{4}$ & $-be_{3}$ \\ 
$e_{3}$ & $e_{3}$ & $-e_{4}$ & $-c$ & $ce_{2}$ \\ 
$e_{4}$ & $e_{4}$ & $be_{3}$ & $-ce_{2}$ & $-bc$%
\end{tabular}%
.\medskip
\end{center}

If $x=x_{1}\cdot 1+x_{2}e_{2}+x_{3}e_{3}+x_{4}e_{4}$$\in $$\mathbb{H}_{%
\mathbb{Q}}\left( -b,-c\right) ,$ we denote by $\boldsymbol{n}\left(
x\right) $ and $\mathbf{t}\left( x\right) ~$\textit{the norm} and \textit{%
the trace} of a quaternion $x$. The \textit{conjugate} of the quaternion $x$
is $\overline{x}=x_{1}\cdot 1-x_{2}e_{2}-x_{3}e_{3}-x_{4}e_{4}$$\in $$%
\mathbb{H}_{\mathbb{Q}}\left( -b,-c\right) $. The trace of a quaternion is 
\begin{equation*}
\mathbf{t}\left( x\right) =x+\overline{x}\in \mathbb{Q}
\end{equation*}%
and the norm has the following expression 
\begin{equation}
\boldsymbol{n}\left( x\right) =x\cdot \overline{x}%
=x_{1}^{2}+bx_{2}^{2}+cx_{3}^{2}+bcx_{4}^{2}.  \tag{3.2.}
\end{equation}%
If, \thinspace for $x\in \mathbb{H}\left( -b,-c\right) $, the relation $%
\mathbf{n}\left( x\right) =0$ implies $x=0$, then the algebra $\mathbb{H}%
\left( -b,-c\right) $ is called a \textit{division} algebra, otherwise the
quaternion algebra is called a \textit{split} algebra.\medskip \newline
\qquad \textbf{Proposition 3.4. }\textit{Let} $b,c$ \textit{be positive
integers and let }$\mathbb{H}_{\mathbb{Q}}\left( -b,-c\right) $\textit{\ be
a generalized quaternion algebra}$.$ If 
\begin{equation*}
\left( b,c\right) \in \left\{ \left( 1,1\right) ,\left( 1,2\right) ,\left(
1,3\right) ,\left( 2,2\right) ,\left( 2,3\right) ,\left( 2,4\right) ,\left(
2,5\right) \right\} ,\newline
\end{equation*}%
\textit{then, the norm map} 
\begin{equation*}
\mathbf{n}:\mathbb{H}_{\mathbb{Q}}\left( -b,-c\right) \rightarrow \mathbb{Q}%
_{+}
\end{equation*}%
\textit{is a surjective function}.\medskip\ \medskip \newline
\qquad \textbf{Proof.} If we consider any positive integer $m$, from
relation (1.1), there are some integers $x_{1},x_{2},x_{3},x_{4}$ such that $%
m=x_{1}^{2}+bx_{2}^{2}+cx_{3}^{2}+bcx_{4}^{2}.$ Thus, there is the
quaternion $x=x_{1}\cdot 1+x_{2}e_{2}+x_{3}e_{3}+x_{4}e_{4}$$\in $$\mathbb{H}%
_{\mathbb{Q}}\left( -b,-c\right) $, such that $m=\mathbf{n}\left( x\right) .$%
\newline
If we consider any positive rational number $m=\frac{m^{^{\prime }}}{l},$
where $m^{^{\prime }},l$$\in $$\mathbb{N},$ $l\neq 0,$ $g.c.d\left(
m^{^{\prime }},l\right) =1$, by using the same argument as above for $%
m^{\prime }l$, there are some rational numbers $x_{1}^{^{\prime }}=\frac{%
x_{1}}{l},x_{2}^{^{\prime }}=\frac{x_{2}}{l},x_{3}^{^{\prime }}=\frac{x_{3}}{%
l},x_{4}^{^{\prime }}=\frac{x_{4}}{l}$ such that\newline
\begin{equation*}
m=\frac{m^{^{\prime }}}{l}=\frac{m^{^{\prime }}l}{l^{2}}=\left( \frac{x_{1}}{%
l}\right) ^{2}+b\left( \frac{x_{2}}{l}\right) ^{2}+c\left( \frac{x_{3}}{l}%
\right) ^{2}+bc\left( \frac{x_{4}}{l}\right) ^{2}=\mathbf{n}\left(
x^{^{\prime }}\right) ,
\end{equation*}%
where $x^{^{\prime }}=x_{1}^{^{\prime }}\cdot 1+x_{2}^{^{\prime
}}e_{2}+x_{3}^{^{\prime }}e_{3}+x_{4}^{^{\prime }}e_{4}$$\in $$\mathbb{H}_{%
\mathbb{Q}}\left( -b,-c\right) $. From here, it results that, $\mathbf{n}$
is a surjective function. $_{{}}\medskip $\medskip \newline

\textbf{Proposition 3.5. }\textit{If} 
\begin{equation*}
\left( b,c\right) \in \left\{ \left( 1,1\right) ,\left( 1,2\right) ,\left(
1,3\right) ,\left( 2,2\right) ,\left( 2,3\right) ,\left( 2,4\right) ,\left(
2,5\right) \right\} ,
\end{equation*}%
\textit{ten the obtained} $7$ \textit{quaternion algebras over the field of
rational numbers,} $\mathbb{H}_{\mathbb{Q}}\left( -b,-c\right) $, \textit{%
are all division algebras}.\medskip

\textbf{Proof.} Since $b>0,$ $c>0,$ it result that the equation 
\begin{equation*}
-bx^{2}-cy^{2}=z_{0}^{2}
\end{equation*}%
does not have solutions $\left( x,y,z\right) $ in $\left( \mathbb{Q}%
\backslash \left\{ \left( 0,0,0\right) \right\} \right) $$\times $$\left( 
\mathbb{Q}\backslash \left\{ \left( 0,0,0\right) \right\} \right) $$\times $$%
\left( \mathbb{Q}\backslash \left\{ \left( 0,0,0\right) \right\} \right) $.$%
_{{}}\medskip \medskip $

Let $K$ be a field. It is known that two quaternion $K$- algebras which
split are isomorphic, but a similar result is not known in the case of two
division quaternion $K-$ algebras. These algebras may be isomorphic or may
be not. We are interested to know how many of these $7$ quaternions $\mathbb{%
Q}$- algebras are isomorphic and which with which is isomorphic.The answer
to this question is given in the following proposition. \newpage 

\textbf{Proposition 3.6.} \newline
\textit{a)} \textit{The quaternion algebras} $\mathbb{H}_{\mathbb{Q}}\left(
-1,-1\right) ,$ $\mathbb{H}_{\mathbb{Q}}\left( -1,-2\right) ,$\newline
$\mathbb{H}_{\mathbb{Q}}\left( -2,-2\right) ,$ $\mathbb{H}_{\mathbb{Q}%
}\left( -2,-3\right) ,$ $\mathbb{H}_{\mathbb{Q}}\left( -2,-4\right) $ 
\textit{are isomorphic}.\newline
\textit{b)} \textit{The quaternion algebras} $\mathbb{H}_{\mathbb{Q}}\left(
-1,-3\right) ,$ $\mathbb{H}_{\mathbb{Q}}\left( -2,-5\right) $ \textit{are
not isomorphic between them and nor with those 5 quaternions algebras from a)%
}.\medskip \newline

\textbf{Proof.} It is known ([Ko]) that if a prime $p$ divides the
discriminant of the algebra $\mathbb{H}_{\mathbb{Q}}\left( -b,-c\right) $,
then $p$ must divides $2bc$. We know that the discriminant $D_{\mathbb{H}_{%
\mathbb{Q}}\left( -b,-c\right) }$ of the algebra $\mathbb{H}_{\mathbb{Q}%
}\left( -b,-c\right) $ is the product of prime ideals of $\mathbb{Z}$ which
ramify in $\mathbb{H}_{\mathbb{Q}}\left( -b,-c\right) $. Since $\left( 
\mathbb{Z},+,\cdot \right) $ is a principal ideal domain, to simplify
writing, an ideal of $\mathbb{Z}$ is identified with its generator, up to
units.\newline
We consider the algebra $\mathbb{H}_{\mathbb{Q}}(-2,-3)$. Let $p$ be a prime
which divides the discriminant of $\mathbb{H}_{\mathbb{Q}}(-2,-3)$,
therefore we have $p$$|$$6.$\newline
If $p=3$, by using the properties of the Hilbert symbol and of the Legendre
symbol we obtain:\newline
$(-2,-3)_{3}=(-1,-3)_{3}(2,-3)_{3}=(-1,-1)_{3}(-1,3)_{3}(2,-1)_{3}(2,3)_{3}=$
\begin{equation*}
=1\cdot \left( \frac{-1}{3}\right) \cdot 1\cdot \left( \frac{2}{3}\right)
=\left( \frac{2}{3}\right) ^{2}=1,
\end{equation*}%
It results that $3$ does not ramify in $H_{\mathbb{Q}}(-2,-3).$\newline
By following computation%
\begin{equation*}
(-2,-3)_{2}=(-1,-3)_{2}(2,-3)_{2}=(-1,-1)_{2}(-1,3)_{2}(2,-1)_{2}(2,3)_{2}=
\end{equation*}%
\begin{equation*}
=\left( -1\right) (-1)^{\frac{3-1}{2}\cdot \frac{-1-1}{2}}1(-1)^{\frac{%
3^{2}-1}{8}}=-1,
\end{equation*}%
we obtain that the only prime which ramifies in $\mathbb{H}_{\mathbb{Q}%
}(-2,-3)$ is $2$. Therefore, the reduced discriminant of the algebra $%
\mathbb{H}_{\mathbb{Q}}(-2,-3)$ is $2.$\newline
In a similar way, we obtain that 
\begin{equation*}
D_{\mathbb{H}_{\mathbb{Q}}\left( -1,-1\right) }=D_{\mathbb{H}_{\mathbb{Q}%
}\left( -1,-2\right) }=D_{\mathbb{H}_{\mathbb{Q}}\left( -2,-2\right) }=D_{%
\mathbb{H}_{\mathbb{Q}}\left( -2,-4\right) }=2.
\end{equation*}%
Now, by using Proposition 3.3, it results that the quaternion algebras $%
\mathbb{H}_{\mathbb{Q}}\left( -1,-1\right) ,$ $\mathbb{H}_{\mathbb{Q}}\left(
-1,-2\right) ,$ $\mathbb{H}_{\mathbb{Q}}\left( -2,-2\right) ,$ $\mathbb{H}_{%
\mathbb{Q}}\left( -2,-3\right) ,$ $\mathbb{H}_{\mathbb{Q}}\left(
-2,-4\right) $ are isomorphic.\newline
b) Similar to a), we calculate the discriminants for the quaternion algebras 
$\mathbb{H}_{\mathbb{Q}}\left( -1,-3\right) ,$ $\mathbb{H}_{\mathbb{Q}%
}\left( -2,-5\right) $ and we obtain that 
\begin{equation*}
D_{\mathbb{H}_{\mathbb{Q}}\left( -1,-3\right) }=3,\ D_{\mathbb{H}_{\mathbb{Q}%
}\left( -2,-5\right) }=5.
\end{equation*}%
Therefore, by using again Proposition 3.3 and a), from above, we obtain that
the quaternion algebras $\mathbb{H}_{\mathbb{Q}}\left( -1,-3\right) ,$ $%
\mathbb{H}_{\mathbb{Q}}\left( -2,-5\right) $ are not isomorphic between them
and nor with those $5$ quaternions algebras from a).$_{{}}\medskip $
\medskip \newline

\textbf{Proposition 3.7. }\textit{Let }$\mathbb{Q}\left( i\right) $ \textit{%
be the quadratic field}, \textit{where} $i^{2}=-1$. \textit{If} 
\begin{equation*}
\left( b,c\right) \in \newline
\left\{ \left( 1,1\right) ,\left( 1,2\right) ,\left( 1,3\right) ,\left(
2,2\right) ,\left( 2,3\right) ,\left( 2,4\right) \right\} ,
\end{equation*}%
\textit{then these} $6$ \textit{quaternion algebras } $\mathbb{H}_{\mathbb{Q}%
\left( i\right) }\left( -b,-c\right) $, \textit{over the quadratic field} $%
\mathbb{Q}\left( i\right) $, \textit{are all split algebras and} $\mathbb{H}%
_{\mathbb{Q}\left( i\right) }\left( -2,-5\right) $ \textit{is a division
algebra}.\medskip

\textbf{Proof.} For the first three cases, we use Proposition 3.1, but we
are not looking for integer solutions of the equation $x^{2}+my^{2}=z^{2}$,
we are looking for solutions over $\mathbb{Z}\left[ i\right] $ of the
equation $-x^{2}-my^{2}=z^{2}.$\newline
\textbf{Case 1. }If $b=1,c=1$, then we remark that $\left(
x_{0},y_{0},z_{0}\right) =\left( 3i,4i,5\right) $ is a solution of the
equation $-x^{2}-y^{2}=z^{2}$ over $\mathbb{Q}\left( i\right) $$\backslash $$%
\left\{ \left( 0,0,0\right) \right\} .$\newline
\textbf{Case 2.} If $b=1,c=2$, then we remark that $\left(
x_{0},y_{0},z_{0}\right) =\left( i,2i,3\right) $ is a solution of the
equation $-x^{2}-2y^{2}=z^{2}$ over $\mathbb{Q}\left( i\right) $$\backslash $%
$\left\{ \left( 0,0,0\right) \right\} .$\newline
\textbf{Case 3.} If $b=1,c=3$, then we remark that $\left(
x_{0},y_{0},z_{0}\right) =\left( i,i,2\right) $ is a solution of the
equation $-x^{2}-3y^{2}=z^{2}$ over $\mathbb{Q}\left( i\right) $$\backslash $%
$\left\{ \left( 0,0,0\right) \right\} .\newline
$\textbf{Case 4.} If $b=2,c=2$, then we remark that $\left(
x_{0},y_{0},z_{0}\right) =\left( i,i,2\right) $ is a solution of the
equation $-2x^{2}-2y^{2}=z^{2}$ over $\mathbb{Q}\left( i\right) $$\backslash 
$$\left\{ \left( 0,0,0\right) \right\} .$\newline
\textbf{Case 5.} If $b=2,c=3$, then we remark that $\left(
x_{0},y_{0},z_{0}\right) =\left( 1,i,1\right) $ is a solution the equation $%
-2x^{2}-3y^{2}=z^{2}$ over $\mathbb{Q}\left( i\right) $$\backslash $$\left\{
\left( 0,0,0\right) \right\} .$ \newline
\textbf{Case 6.} If $b=2,c=4$, then we remark that $\left(
x_{0},y_{0},z_{0}\right) =\left( 2i,1,2\right) $ is a solution of the
equation $-2x^{2}-4y^{2}=z^{2}$over $\mathbb{Q}\left( i\right) $$\backslash $%
$\left\{ \left( 0,0,0\right) \right\} .$\newline
Therefore, from Proposition 3.2, it results that the quaternions algebras $%
\mathbb{H}_{\mathbb{Q}\left( i\right) }\left( -1,-1\right),\newline
\mathbb{H}_{\mathbb{Q}\left( i\right) }\left( -1,-2\right) ,$ $\mathbb{H}_{%
\mathbb{Q}\left( i\right) }\left( -1,-3\right) ,$ $\mathbb{H}_{\mathbb{Q}%
\left( i\right) }\left( -2,-2\right) ,$ $\mathbb{H}_{\mathbb{Q}\left(
i\right) }\left( -2,-3\right) ,$ $\mathbb{H}_{\mathbb{Q}\left( i\right)
}\left( -2,-4\right) $ are split algebras.\newline
\textbf{Case 7.} If $b=2,c=5$, we use the same idea as in Proposition 3.6
and we determine the primes which ramify in $\mathbb{H}_{\mathbb{Q}\left(
i\right) }\left( -2,-5\right) $. For the same purpose, we also can use the
computer algebra system Magma ([Mag]) and we find that the reduced
discriminat of the quaternion algebra $\mathbb{H}_{\mathbb{Q}\left( i\right)
}\left( -b,-c\right) $ is $D_{\mathbb{H}_{\mathbb{Q}\left( i\right) }\left(
-2,-5\right) }=5\mathbb{Z}\left[ i\right] $, therefore the quaternion
algebra $\mathbb{H}_{\mathbb{Q}\left( i\right) }\left( -2,-5\right) $ is a
division algebra. $_{{}}\medskip $

\begin{equation*}
\end{equation*}

\bigskip \textbf{4. Some aspects regarding the arithmetic of Integer
Quaternions}%
\begin{equation*}
\end{equation*}

As we remarked above, the algebras $\mathbb{H}_{\mathbb{Q}}\left(
-1,-1\right) ,$ $\mathbb{H}_{\mathbb{Q}}\left( -1,-2\right) ,$\newline
$\mathbb{H}_{\mathbb{Q}}\left( -2,-2\right) ,$ $\mathbb{H}_{\mathbb{Q}%
}\left( -2,-3\right) ,$ $\mathbb{H}_{\mathbb{Q}}\left( -2,-4\right) $, $%
\mathbb{H}_{\mathbb{Q}}\left( -1,-3\right) ,$ $\mathbb{H}_{\mathbb{Q}}\left(
-2,-5\right) $ are division algebras and have a surjective norm form. We
denote with $\mathcal{H}$ the set of these algebras, namely $\mathcal{H}=\{%
\mathbb{H}_{\mathbb{Q}}\left( -1,-1\right) ,\mathbb{H}_{\mathbb{Q}}\left(
-1,-2\right) ,\newline
\mathbb{H}_{\mathbb{Q}}\left( -2,-2\right) ,\mathbb{H}_{\mathbb{Q}}\left(
-2,-3\right) ,\mathbb{H}_{\mathbb{Q}}\left( -2,-4\right) ,\mathbb{H}_{%
\mathbb{Q}}\left( -1,-3\right) ,\mathbb{H}_{\mathbb{Q}}\left( -2,-5\right) \}
$. First five of these algebras are isomophic. In $\mathbb{H}_{\mathbb{Q}%
}\left( -1,-1\right) $ is defined the following sets 
\begin{equation*}
\mathbb{H}\left( \mathbb{Z}\right) \text{=}\{x\in \mathbb{H}_{\mathbb{Q}%
}\left( -1,-1\right) ,x\text{=}x_{1}\cdot
1+x_{2}e_{2}+x_{3}e_{3}+x_{4}e_{4},x_{1},x_{2},x_{3},x_{4}\in \mathbb{Z}\},
\end{equation*}%
called the set of \textit{Lipschitz integers} and%
\begin{eqnarray*}
\mathbb{H}\left( \mathbb{Z}\left[ \frac{1}{2}\right] \right)  &\text{=}%
&\{x\in \mathbb{H}_{\mathbb{Q}}\left( -1,-1\right) ,x\text{=}x_{1}\cdot
1+x_{2}e_{2}+x_{3}e_{3}+x_{4}e_{4},x_{1},x_{2},x_{3},x_{4}\in \mathbb{Z}\},
\\
&&x_{1},x_{2},x_{3},x_{4}\text{ are all in }\mathbb{Z}\text{ or all in }%
\mathbb{Z}+\frac{1}{2},
\end{eqnarray*}%
called the set of \textit{Hurwitz integers} (see [Co, Sm; 03], p.55).

In a similar way, for an algebra $A\in \mathcal{H}$, we will define the set 
\begin{equation*}
A\left( \mathbb{Z}\right) =\{x\in A,x\text{=}x_{1}\cdot
1+x_{2}e_{2}+x_{3}e_{3}+x_{4}e_{4},x_{1},x_{2},x_{3},x_{4}\in \mathbb{Z}\}.
\end{equation*}%
We call the elements of this set \textit{the integers of the algebra} $A\in 
\mathcal{H}$.\medskip 

\textbf{Proposition 4.1.} \textit{For an algebra} $A\in \mathcal{H}$\textit{%
, for} $x,y\in A\left( \mathbb{Z}\right) $\textit{, there are} $\gamma $ 
\textit{and} $\theta $ \textit{such that} 
\begin{equation*}
n\left( y\right) x=\gamma y+\theta n\left( y\right) ,
\end{equation*}%
\textit{with} $n\left( \theta \right) <n\left( y\right) ,\gamma ,\theta \in
A\left( \mathbb{Z}\right) .\medskip $

\textbf{Proof.} Since the norm $n$ is surjective, there is a quaternion $%
\theta \in A\left( \mathbb{Z}\right) $, such that $n\left( \theta \right)
<n\left( y\right) $. Therefore, there is a quaternion $\delta \in A$ such
that $\theta =x-\delta y$, with 
\begin{equation*}
\delta =\left( x-\theta \right) y^{-1}=\frac{(x-\theta )\overline{y}}{n(y)}
\end{equation*}%
We denote $\gamma =(x-\theta )\overline{y}\in A\left( \mathbb{Z}\right) $.
It results 
\begin{equation*}
x=\delta y+\theta =\frac{\gamma }{n\left( y\right) }y+\theta ,
\end{equation*}%
therefore $n\left( y\right) x=\gamma y+n\left( y\right) \theta $. $%
_{{}}\medskip $

\textbf{Proposition 4.2.} ( [Da, Sa, Va; 03], Lemma 2.6.3). \textit{Let} $%
x,y\in \mathbb{H}\left( \mathbb{Z}\right) $\textit{, with} $n\left( y\right) 
$ \textit{an odd number. Therefore, there are} $\gamma $ \textit{and} $%
\theta \in \mathbb{H}\left( \mathbb{Z}\right) $ \textit{such that} 
\begin{equation*}
x=\gamma y+\theta ,\text{ \textit{with} }n\left( \theta \right) <n\left(
y\right) .
\end{equation*}%
$_{{}}$

\textbf{Remark 4.3.} Proposition 4.2 was proved without using in the proof
the surjectivity of the norm. Proposition 4.2 is a particular case of the
Proposition 4.1, since if $n\left( y\right) $ is arbitrary, then the result
obtained in Proposition 4.2\ is a;so true in $\mathbb{H}\left( \mathbb{Z}%
\left[ \frac{1}{2}\right] \right) $.\medskip

\textbf{Definition 4.4.} Let $q\in \mathbb{H}\left( \mathbb{Z}\right) $. If $%
\mathbf{n}\left( q\right) $ is an odd integer, then the quaternion $q$ is
called an \textit{odd} quaternion.\medskip\ 

\textbf{Definition 4.5. }([Gu; 13], Definition 2 and Theorem 1) We consider $%
\varphi $ an odd quaternion from $\mathbb{H}\left( \mathbb{Z}\right) $ and $%
x,y\in \mathbb{H}\left( \mathbb{Z}\right) $ two arbitrary quaternions. If
there is the quaternion $\theta \in \mathbb{H}\left( \mathbb{Z}\right) $
such that $x-y=\theta \varphi $, therefore the quaternions $x$ and $y$ are
called \textit{right congruent modulo} $\varphi $ and we will denote this $%
x\equiv _{r}y$ mod $\varphi $. We obtain an equivalence relation and the
quotient set is 
\begin{equation*}
\mathbb{H}\left( \mathbb{Z}\right) _{\varphi }=\{x\text{ mod }\varphi ,x\in 
\mathbb{H}\left( \mathbb{Z}\right) \}.
\end{equation*}%
The set $\mathbb{H}\left( \mathbb{Z}\right) _{\varphi }$ has $n^{2}\left(
\varphi \right) $ elements and, from the above propositions, we have $%
n\left( y\right) <n\left( \varphi \right) $.\medskip

\textbf{Definition 4.6}. Let $A\in \mathcal{H}$, a quaternion $q\in A$ is 
\textit{a unit} in $A$ if $n\left( q\right) =1$. A quaternion $q\in A\left( 
\mathbb{Z}\right) ~$is a \textit{prime quaternion} if $q=\alpha \beta
,\alpha ,\beta \in A$ implies $\alpha $ or $\beta $ is a unit
element.\medskip

\textbf{Remark 4.7}. ([Gu; 13] , Definition 5 and Remark 2)

1) We will consider the element $v\in \mathbb{H}\left( \mathbb{Z}\left[ 
\frac{1}{2}\right] \right) ,v=\frac{1}{2}\left( 1+e_{2}+e_{3}+e_{4}\right) $
and we consider the set $\mathcal{A}=\{q\in \mathbb{H}\left( \mathbb{Z}\left[
\frac{1}{2}\right] \right) ,q=\alpha +\beta v,\alpha ,\beta \in \mathbb{Z}\}$%
. We remark that multiplication is commutative on $\mathcal{A}$, therefore
we obtain that this set is a commutative ring. It results that we can define
the same equivalence relation as above, obtaining the quotient ring $%
\mathcal{A}_{\varphi }$, with $\varphi $ a prime quaternion in $\mathcal{A}$%
.\medskip

\textbf{Proposition 4.8. }\textit{A quaternion} $q\in A\left( \mathbb{Z}%
\right) ~$\textit{is a prime quaternion if and only if its norm} $\mathbf{n}%
\left( q\right) $ \textit{is a prime number}.\medskip

\textbf{Proof.} We supose that $q$ is a prime quaternion and $\mathbf{n}%
\left( q\right) =m$, such that $m$ is not a prime number.Therefore, $m=rs$,
with $r\neq 1$ and $s\neq 1$. Since the norm $\mathbf{n}$ is surjective, it
results that there are the quaternions $q_{1},q_{2}\in A\left( \mathbb{Z}%
\right) $ such that $\mathbf{n}\left( q_{1}\right) =r$ and $\mathbf{n}\left(
q_{2}\right) =s$. We obtain that $q_{1}$ and $q_{2}$ are not unit
quaternions, thus $q$ is not a prime quaternion, which it is false. It
results that $\mathbf{n}\left( q\right) $ is a prime number.

Conversely, if $\mathbf{n}\left( q\right) =p$ is a prime number and $q=$ $%
q_{1}q_{2}$, we obtain that $p=\mathbf{n}(q)=\mathbf{n}($ $q_{1}q_{2})=%
\mathbf{n}\left( q_{1}\right) \mathbf{n}\left( q_{2}\right) $. Therefore, $%
\mathbf{n}\left( q_{1}\right) =1$ or \ $\mathbf{n}\left( q_{2}\right) =1$.
It results that $q$ is a prime quaternion.$_{{}}\medskip $

\textbf{Remark 4.9.} The above proposition generalized Corollary 2.6.10 from
[Da, Sa, Va; 03]. In that statement the above result was proved for the
algebra $\mathbb{H}_{\mathbb{Q}}\left( -1,-1\right) $ without using the
surjectivity of the norm.

2) ([Gu; 13], Theorem 2) If $q\in \mathcal{A}$, $q==\alpha +\beta v,\alpha
,\beta \in \mathbb{Z}$, with $\gcd \{\alpha ,\beta \}=1$, therefore $\mathbb{%
Z}_{n\left( q\right) }$ and $\mathcal{A}_{\varphi }$ are isomorphic as
rings. \ Therefore, if $q$ is a prime quaternion, that means its norm $%
\mathbf{n}\left( q\right) $ is a prime number, we obtain than $\mathcal{A}%
_{\varphi }$ is a commutative field isomorphic to $\mathbb{Z}_{n\left(
\varphi \right) }$. We remark that 
\begin{eqnarray*}
n\left( \varphi \right) &=&\varphi \overline{\varphi }=\left( \alpha +\beta
v\right) \left( \alpha +\beta \overline{v}\right) = \\
&=&\alpha ^{2}+\alpha \beta \overline{v}+\alpha \beta v+\beta ^{2}v\overline{%
v}= \\
&=&\alpha ^{2}+\alpha \beta +\beta ^{2}\text{.}
\end{eqnarray*}

\textbf{Example 4.10.} 1) Let \thinspace $q=-1+2v$. We obtain that $n\left(
q\right) =3$, therefore $\mathbb{Z}_{3}\simeq \mathcal{A}_{-1+2v}=\{0,1,v\}$%
, since we have

$2=(-1+2v)(-v)+v$, with the norm $n\left( v\right) =1<2.$

2) Let $q=1+2v$. We obtain that $n\left( q\right) =7,$ therefore $\mathbb{Z}%
_{7}\simeq \mathcal{A}_{1+2v}=$\newline
$=\{0,1,v-1,v,1+v,2v-1,2v\}$.

We have:

$2=\left( 1+2v\right) \left( 1-v\right) +v-1$, with the norm $n\left(
v-1\right) =1<7$.

$3=\left( 1+2v\right) \left( 1-v\right) +v$, with the norm $n\left( v\right)
=1<7$.

$4=\left( 1+2v\right) \left( 1-v\right) +1+v$, with the norm $n\left(
1+v\right) =\frac{12}{4}=3<7;$

$5=\left( 1+2v\right) \left( 2-2v\right) +2v-1$, with the norm $n\left(
2v-1\right) =3<7;$

$6=$ $\left( 1+2v\right) \left( 2-2v\right) +2v$, with the norm $n\left(
2v\right) =4<7.\medskip $

\textbf{Remark 4.11.} The isomorphism beween $\mathbb{Z}_{n\left( q\right) }$
and $\mathcal{A}_{\varphi }$ has direct applications in Coding Theory. For
other details, see for example [GF; 10], [Hu; 94], [Ri; 95].

\begin{equation*}
\end{equation*}

\bigskip \textbf{5. Quaternions with permutations} 
\begin{equation*}
\end{equation*}

For $b=1,c=2,$ we consider the quaternion algebra $\mathbb{H}_{\mathbb{Q}%
\left( i\right) }\left( -1,-2\right) $. From Proposition 3.7, this algebra
is a split algebra.\newline

Let $S_{4}$ be the set of all permutations of degree $4$ on the set $\left\{
n,n+1,n+2,n+3\right\} $, where $n$ is a positive integer. We define a new
element, namely, the $\sigma -$ \textit{permutated }$n$\textit{th
Fibonacci-Hurtwitz quaternion} 
\begin{equation*}
F_{n,\sigma }=\frac{1}{2}f_{\sigma \left( n\right) }+\frac{1}{2}f_{\sigma
\left( n+1\right) }e_{2}+\frac{1}{2}f_{\sigma \left( n+2\right) }e_{2}+\frac{%
1}{2}f_{\sigma \left( n+3\right) }e_{4}.\newline
\end{equation*}%
\medskip \newline

\textbf{Proposition 5.1.} \textit{Let } $\mathbb{Q}\left( i\right) $ \textit{%
be the quadratic field} \textit{where} $i^{2}=-1$. \textit{For the
quaternion algebra} $\mathbb{H}_{\mathbb{Q}\left( i\right) }\left(
-1,-2\right) $ \textit{and for each positive integer} $n$ \textit{let} $%
F_{n,\sigma }$ \textit{be the} $\sigma -$ \textit{permutated} $n$ \textit{th
Fibonacci-Hurtwitz quaternion}. \textit{Therefore, for} $\sigma =$$%
\begin{pmatrix}
n & n+1 & n+2 & n+3 \\ 
n+3 & n & n+1 & n+2%
\end{pmatrix}%
$ \textit{or for} 
\begin{equation*}
\sigma \in \left\{ 
\begin{array}{c}
\begin{pmatrix}
n & n+1 & n+2 & n+3 \\ 
n & n+3 & n+1 & n+2%
\end{pmatrix}%
,%
\begin{pmatrix}
n & n+1 & n+2 & n+3 \\ 
n+3 & n & n+2 & n+1%
\end{pmatrix}%
, \\ 
\newline
\begin{pmatrix}
n & n+1 & n+2 & n+3 \\ 
n & n+3 & n+2 & n+1%
\end{pmatrix}%
\end{array}%
\right\}
\end{equation*}%
\textit{the norm of a such quaternion is} 
\begin{equation*}
\mathbf{n}\left( F_{n,\sigma }\right) =f_{2n+3}.
\end{equation*}%
\medskip \newline
\textbf{Proof.} For $\sigma =$$%
\begin{pmatrix}
n & n+1 & n+2 & n+3 \\ 
n+3 & n & n+1 & n+2%
\end{pmatrix}%
$, we have 
\begin{equation*}
\mathbf{n}\left( F_{n,\sigma }\right) =\frac{1}{4}%
(f_{n+3}^{2}+f_{n}^{2}+2f_{n+1}^{2}+3f_{n+2}^{2}).
\end{equation*}%
Using Proposition 2.1 i) and ii), we obtain 
\begin{equation*}
\mathbf{n}\left( F_{n,\sigma }\right) =\frac{1}{4}%
(4f_{n+1}^{2}+4f_{n+2}^{2})=f_{n+1}^{2}+f_{n+2}^{2}=f_{2n+3}.
\end{equation*}%
In the same way, for\newline
\begin{equation*}
\sigma \in \left\{ 
\begin{array}{c}
\begin{pmatrix}
n & n+1 & n+2 & n+3 \\ 
n & n+3 & n+1 & n+2%
\end{pmatrix}%
,%
\begin{pmatrix}
n & n+1 & n+2 & n+3 \\ 
n+3 & n & n+2 & n+1%
\end{pmatrix}%
, \\ 
\newline
\begin{pmatrix}
n & n+1 & n+2 & n+3 \\ 
n & n+3 & n+2 & n+1%
\end{pmatrix}%
\end{array}%
\right\} ,
\end{equation*}%
we obtain that $\mathbf{n}\left( F_{n,\sigma }\right) =f_{2n+3}$. $%
_{{}}\medskip $ \medskip \newline

\textbf{Corollary 5.2.} \textit{We consider } $\mathbb{Q}\left( i\right) $%
\textit{\ the quadratic field} \textit{where} $i^{2}=-1$ \textit{and the
quaternion algebra} $\mathbb{H}_{\mathbb{Q}\left( i\right) }\left(
-1,-2\right) .$ \textit{Then, for} $\sigma =$ $%
\begin{pmatrix}
n & n+1 & n+2 & n+3 \\ 
n+3 & n & n+1 & n+2%
\end{pmatrix}%
$ \textit{or for} 
\begin{equation*}
\sigma \in \left\{ 
\begin{array}{c}
\begin{pmatrix}
n & n+1 & n+2 & n+3 \\ 
n & n+3 & n+1 & n+2%
\end{pmatrix}%
,%
\begin{pmatrix}
n & n+1 & n+2 & n+3 \\ 
n+3 & n & n+2 & n+1%
\end{pmatrix}%
,\newline
\\ 
\begin{pmatrix}
n & n+1 & n+2 & n+3 \\ 
n & n+3 & n+2 & n+1%
\end{pmatrix}%
\end{array}%
\right\}
\end{equation*}%
\newline
\textit{the following statements we true:\newline
}i) \textit{All }$\sigma -$ \textit{permutated} $n$ \textit{th Fibonacci
quaternions} $F_{n,\sigma }$ \textit{are invertible};\newline
ii) \textit{There are not units}$\sigma -$ \textit{permutated} $n$ \textit{%
th Fibonacci quaternions} $F_{n,\sigma }$.\medskip \newline
\medskip \newline

\textbf{Proof.} i) Applying Proposition 5.1, for $\sigma =$ $%
\begin{pmatrix}
n & n+1 & n+2 & n+3 \\ 
n+3 & n & n+1 & n+2%
\end{pmatrix}%
$ or for 
\begin{equation*}
\sigma \in \left\{ 
\begin{array}{c}
\begin{pmatrix}
n & n+1 & n+2 & n+3 \\ 
n & n+3 & n+1 & n+2%
\end{pmatrix}%
,%
\begin{pmatrix}
n & n+1 & n+2 & n+3 \\ 
n+3 & n & n+2 & n+1%
\end{pmatrix}%
, \\ 
\begin{pmatrix}
n & n+1 & n+2 & n+3 \\ 
n & n+3 & n+2 & n+1%
\end{pmatrix}%
\end{array}%
\right\} ,
\end{equation*}%
we have $\mathbf{n}\left( F_{n,\sigma }\right) =f_{2n+3}\neq 0,$ $\left(
\forall \right) $$n$$\in $$\mathbb{N}$, therefore $F_{n,\sigma }$ are
invertible \ for all $n$$\in $$\mathbb{N}.$\newline
ii) The only Fibonacci numbers equals with $1$ are $f_{1}$ and $f_{2},$
then, by using Proposition 5.1, we have $\mathbf{n}\left( F_{n,\sigma
}\right) =f_{2n+3}\neq 1$, for all $n$$\in $$\mathbb{N}$. It results that $%
F_{n,\sigma }$ is not a unit quaternion, for all $n$$\in $$\mathbb{N}$.$%
_{{}}\medskip $

If $n,k$ are two positive integers, from [Ra; 1917], p. 171, we get tha
there are $x_{1},x_{2},x_{3},x_{4},y_{1},y_{2},y_{3},y_{4}$ $\in $$\mathbb{Z}
$ such that $n=x_{1}^{2}+x_{2}^{2}+2x_{3}^{2}+2x_{4}^{2}$ and $%
k=y_{1}^{2}+y_{2}^{2}+2y_{3}^{2}+2y_{4}^{2}.$ We want find $%
u_{1},u_{2},u_{3},u_{4}$$\in $$\mathbb{Z}$ such that $n\cdot
k=u_{1}^{2}+u_{2}^{2}+2u_{3}^{2}+2u_{4}^{2}.$\medskip \newline

\textbf{Proposition 5.3.} \textit{Let} $n,k$ \textit{be two positive
integers,} $n=x_{1}^{2}+x_{2}^{2}+2x_{3}^{2}+2x_{4}^{2}$, $%
k=y_{1}^{2}+y_{2}^{2}+2y_{3}^{2}+2y_{4}^{2}$, \textit{with} $%
x_{1},x_{2},x_{3},x_{4},y_{1},y_{2},y_{3},y_{4}\in $ $\mathbb{Z}.$ \textit{%
Then} 
\begin{equation*}
n\cdot k=\left( x_{1}y_{1}-x_{2}y_{2}-2x_{3}y_{3}-2x_{4}y_{4}\right)
^{2}+\left( x_{1}y_{2}+x_{2}y_{1}+2x_{3}y_{4}-2x_{4}y_{3}\right) ^{2}+
\end{equation*}%
\begin{equation*}
+2\left( x_{1}y_{3}+x_{3}y_{1}-x_{2}y_{4}+x_{4}y_{2}\right) ^{2}+2\left(
x_{1}y_{4}+x_{4}y_{1}+x_{2}y_{3}-x_{3}y_{2}\right) ^{2}.
\end{equation*}

\textbf{Proof.} We consider the quaternion algebra $\mathbb{H}_{\mathbb{Q}%
}\left( -1,-2\right) $ and $x=x_{1}\cdot 1+x_{2}e_{2}+x_{3}e_{3}+x_{4}e_{4}$%
, $y=y_{1}\cdot 1+y_{2}e_{2}+y_{3}e_{3}+y_{4}e_{4}$$\in $$\mathbb{H}_{%
\mathbb{Q}}\left( -1,-2\right) ,$ where $%
x_{1},x_{2},x_{3},x_{4},y_{1},y_{2},y_{3},y_{4}$$\in $$\mathbb{Z}$. We have:%
\newline
\begin{equation*}
x\cdot y=\left( x_{1}y_{1}-x_{2}y_{2}-2x_{3}y_{3}-2x_{4}y_{4}\right) \cdot
1+\left( x_{1}y_{2}+x_{2}y_{1}+2x_{3}y_{4}-2x_{4}y_{3}\right) e_{1}+
\end{equation*}%
\begin{equation*}
+\left( x_{1}y_{3}+x_{3}y_{1}-x_{2}y_{4}+x_{4}y_{2}\right) e_{2}+\left(
x_{1}y_{4}+x_{4}y_{1}+x_{2}y_{3}-x_{3}y_{2}\right) e_{3}.
\end{equation*}%
By using the fact that the norm is multiplicative, that means $\boldsymbol{n}%
\left( xy\right) =\boldsymbol{n}\left( x\right) \boldsymbol{n}\left(
y\right) $, we obtain 
\begin{equation*}
n\cdot k=\left( x_{1}y_{1}-x_{2}y_{2}-2x_{3}y_{3}-2x_{4}y_{4}\right)
^{2}+\left( x_{1}y_{2}+x_{2}y_{1}+2x_{3}y_{4}-2x_{4}y_{3}\right) ^{2}+
\end{equation*}%
\begin{equation*}
+2\left( x_{1}y_{3}+x_{3}y_{1}-x_{2}y_{4}+x_{4}y_{2}\right) ^{2}+2\left(
x_{1}y_{4}+x_{4}y_{1}+x_{2}y_{3}-x_{3}y_{2}\right) ^{2}.
\end{equation*}%
Since the numbers $x_{1},x_{2},x_{3},x_{4}$, obtained in the expresion of $n$%
, are not unique and the numbers $y_{1},y_{2},y_{3},y_{4}$, obtained in the
expresion of $k$, are not unique, the representation of $nk$ under the form $%
n\cdot k=u_{1}^{2}+u_{2}^{2}+2u_{3}^{2}+2u_{4}^{2}$ is not unique. $_{{}}$
\medskip \newline

\textbf{Corollary 5.4.} \textit{Let} $n,k$ \textit{be two positive integers,}
$n=x_{1}^{2}+x_{2}^{2}+2x_{3}^{2}+2x_{4}^{2},$ $%
k=y_{1}^{2}+y_{2}^{2}+2y_{3}^{2}+2y_{4}^{2},$ \textit{with} $%
x_{1},x_{2},x_{3},x_{4},y_{1},y_{2},y_{3},y_{4}\in $ $\mathbb{Z}$. \textit{%
Then we have} 
\begin{equation*}
n\cdot k=\left( x_{1}y_{1}+x_{2}y_{2}+2x_{3}y_{3}+2x_{4}y_{4}\right)
^{2}+\left( x_{1}y_{2}-x_{2}y_{1}+2x_{3}y_{4}-2x_{4}y_{3}\right) ^{2}+
\end{equation*}%
\begin{equation*}
+2\left( x_{1}y_{3}-x_{3}y_{1}-x_{2}y_{4}+x_{4}y_{2}\right) ^{2}+2\left(
x_{1}y_{4}-x_{4}y_{1}+x_{2}y_{3}-x_{3}y_{2}\right) ^{2}.
\end{equation*}%
\textbf{Proof.} By straightforward calculation. $_{{}}$ \medskip \newline

\textbf{Remark 5.5.} The expresion of $nk$ obtained in Corollary 5.4
corresponds to the norm of the product of the following quaternions 
\begin{equation*}
x^{^{\prime }}=x_{1}\cdot 1-x_{2}e_{2}-x_{3}e_{3}+x_{4}e_{4}\in \mathbb{H}_{%
\mathbb{Q}}\left( -1,-2\right) ,
\end{equation*}%
and%
\begin{equation*}
y^{^{\prime }}=y_{1}\cdot 1+y_{2}e_{2}+y_{3}e_{3}-y_{4}e_{4}\in \mathbb{H}_{%
\mathbb{Q}}\left( -1,-2\right) ,
\end{equation*}
where $x_{1},x_{2},x_{3},x_{4},y_{1},y_{2},y_{3},y_{4}$$\in $$\mathbb{Z}$%
.\medskip

From the above results, a question arise: if there is a kind of Fibonacci
elements, possible similar with $\sigma -$ permutated Fibonacci -Hurtwitz
quaternions, such that the product of such elements is also an element of
this type. In this regard, in the following, we provide a nice property
involving the product of two $\sigma -$ permutated Fibonacci -Hurtwitz
quaternions. First of all, we give another property of the Fibonacci
numbers, which we will use later in our proof.\medskip

\textbf{Proposition 5.6.} \textit{Let} $(f_{n})_{n\geq 0}$ \textit{be the
Fibonacci sequence. Then, the following property holds} 
\begin{equation*}
f_{n}f_{l}+f_{n+3}f_{l+3}=2f_{n+l+3},\text{for all}~n,l\in \mathbb{N}.
\end{equation*}

\textbf{Proof.} By using Binet's formula for Fibonacci sequence and Vi\`{e}%
te relations for $\alpha $ and $\beta $, we have 
\begin{equation*}
f_{n}f_{l}+f_{n+3}f_{l+3}=\frac{\alpha ^{n}-\beta ^{n}}{\alpha -\beta }\cdot 
\frac{\alpha ^{l}-\beta ^{l}}{\alpha -\beta }+\frac{\alpha ^{n+3}-\beta
^{n+3}}{\alpha -\beta }\cdot \frac{\alpha ^{l+3}-\beta ^{l+3}}{\alpha -\beta 
}=
\end{equation*}%
\begin{equation*}
=\frac{\alpha ^{n+l+6}+\beta ^{n+l+6}-\alpha ^{n+3}\beta ^{l+3}-\alpha
^{l+3}\beta ^{n+3}+\alpha ^{n+l}+\beta ^{n+l}-\alpha ^{n}\beta ^{l}-\alpha
^{l}\beta ^{n}}{\left( \alpha -\beta \right) ^{2}}=
\end{equation*}%
\begin{equation*}
=\frac{\alpha ^{n+l+3}\left( \alpha ^{3}+\alpha ^{-3}\right) +\beta
^{n+l+3}\left( \beta ^{3}+\beta ^{-3}\right) -\alpha ^{n}\beta ^{l}\left(
\alpha ^{3}\beta ^{3}+1\right) -\alpha ^{l}\beta ^{n}\left( \alpha ^{3}\beta
^{3}+1\right) }{\left( \alpha -\beta \right) ^{2}}=
\end{equation*}%
\begin{equation*}
=\frac{\alpha ^{n+l+3}\left( \alpha ^{3}-\beta ^{3}\right) +\beta
^{n+l+3}\left( \beta ^{3}-\alpha ^{3}\right) }{\left( \alpha -\beta \right)
^{2}}=\frac{\left( \alpha ^{3}-\beta ^{3}\right) \left( \alpha
^{n+l+3}-\beta ^{n+l+3}\right) }{\left( \alpha -\beta \right) ^{2}}=
\end{equation*}%
\begin{equation*}
=\frac{\left( \alpha ^{2}+\alpha \beta +\beta ^{2}\right) \left( \alpha
^{n+l+3}-\beta ^{n+l+3}\right) }{\alpha -\beta }=2f_{n+l+3.}
\end{equation*}%
\smallskip \newline

\textbf{Proposition 5.7.} \textit{Let} $(f_{n})_{n\geq 0}$ \textit{be the
Fibonacci sequence, let} $\sigma =$ $%
\begin{pmatrix}
n & n+1 & n+2 & n+3 \\ 
n+3 & n & n+1 & n+2%
\end{pmatrix}%
$ \textit{and, for each positive integers} $n<l$, \textit{let} $F_{n,\sigma
}^{^{\prime }}$ $and$ $F_{l,\sigma }^{^{\prime \prime }}$ \textit{be the} $%
\sigma -$ \textit{special permutated} \textit{Fibonacci -Hurtwitz quaternions%
}: $F_{n,\sigma }^{^{\prime }}=\frac{1}{2}f_{\sigma \left( n\right) }-\frac{1%
}{2}f_{\sigma \left( n+1\right) }e_{2}-\frac{1}{2}f_{\sigma \left(
n+2\right) }e_{3}+\frac{1}{2}f_{\sigma \left( n+3\right) }e_{4}$, $%
F_{l,\sigma }^{^{\prime \prime }}=\frac{1}{2}f_{\sigma \left( l\right) }+%
\frac{1}{2}f_{\sigma \left( l+1\right) }e_{2}+\frac{1}{2}f_{\sigma \left(
l+2\right) }e_{3}-\frac{1}{2}f_{\sigma \left( l+3\right) }e_{4}$. \textit{Let%
} $k=l-n$. \textit{Therefore, the expression} 
\begin{equation*}
F_{n,\sigma }^{^{\prime }}F_{l,\sigma }^{^{\prime \prime }}-\frac{\mathbf{t}%
\left( F_{n,\sigma }^{^{\prime }}F_{l,\sigma }^{^{\prime \prime }}\right) }{2%
}
\end{equation*}%
\textit{depends only of} $k$.\medskip

\textbf{Proof.} For $F_{n,\sigma }^{^{\prime }}=\frac{1}{2}f_{n+3}-\frac{1}{2%
}f_{n}e_{2}-\frac{1}{2}f_{n+1}e_{3}+\frac{1}{2}f_{n+2}e_{4}$ and\newline
$F_{l,\sigma }^{^{\prime \prime }}=\frac{1}{2}f_{l+3}+\frac{1}{2}f_{l}e_{2}+%
\frac{1}{2}f_{l+1}e_{3}-\frac{1}{2}f_{l+2}e_{4}$ we apply Corollary 4.4 and
Remark 4.5 and we have 
\begin{equation*}
F_{n,\sigma }^{^{\prime }}F_{l,\sigma }^{^{\prime \prime }}=\frac{1}{4}%
\left( f_{n+3}f_{l+3}+f_{n}f_{l}+2f_{n+1}f_{l+1}+2f_{n+2}f_{l+2}\right) +
\end{equation*}%
\begin{equation*}
+\frac{1}{4}\left(
f_{n+3}f_{l}-f_{n}f_{l+3}+2f_{n+1}f_{l+2}-2f_{n+2}f_{l+1}\right) e_{2}+
\end{equation*}%
\begin{equation*}
+\frac{1}{4}\left(
f_{n+3}f_{l+1}-f_{n+1}f_{l+3}-f_{n}f_{l+2}+f_{n+2}f_{l}\right) e_{3}+
\end{equation*}%
\begin{equation*}
+\frac{1}{4}\left(
-f_{n+3}f_{l+2}+f_{n+2}f_{l+3}-f_{n}f_{l+1}+f_{n+1}f_{l}\right) e_{4}.
\end{equation*}%
From the properties of Fibonacci numbers and Remark 5.6., we get 
\begin{equation*}
F_{n,\sigma }^{^{\prime }}F_{l,\sigma }^{^{\prime \prime }}=\frac{1}{4}%
\left( 2f_{n+l+3}+2f_{n+l+3}\right) +\frac{1}{4}\left[ \left( -1\right)
^{n-3}f_{3}f_{l-n}+2\left( -1\right) ^{n}f_{1}f_{l-n+2}\right] e_{2}+
\end{equation*}%
\begin{equation*}
+\frac{1}{4}\left[ \left( -1\right) ^{n+1}f_{2}f_{l-n}+\left( -1\right)
^{n}f_{1}f_{l-n}\right] e_{3}+\frac{1}{4}\left[ \left( -1\right)
^{n+2}f_{1}f_{l-n}+\left( -1\right) ^{n}f_{1}f_{l-n}\right] e_{4}
\end{equation*}%
\begin{equation*}
=f_{n+l+3}+\frac{1}{2}\left( -1\right) ^{n}f_{l-n+1}e_{2}+\frac{1}{2}\left(
-1\right) ^{n}f_{l-n}e_{4}.
\end{equation*}%
It results that $\mathbf{t}(F_{n,\sigma }^{^{\prime }}F_{l,\sigma
}^{^{\prime \prime }})=2f_{n+l+3}$. Therefore, we obtain that 
\begin{equation*}
F_{n,\sigma }^{^{\prime }}F_{l,\sigma }^{^{\prime \prime }}-\frac{\mathbf{t}%
\left( F_{n,\sigma }^{^{\prime }}F_{l,\sigma }^{^{\prime \prime }}\right) }{2%
}=\frac{1}{2}\left( -1\right) ^{n}f_{k+1}e_{2}+\frac{1}{2}\left( -1\right)
^{n}f_{k}e_{4}.
\end{equation*}%
$_{{}}$

\begin{equation*}
\end{equation*}

\textbf{Conclusions. }In this paper we have found an interesting connection
between quaternion algebras and some quaternary quadratic forms, arising
from some results obtained by Ramanujan in \textbf{[}Ra; 1917\textbf{],} we
have define a monoid structure over a finite set on which we will prove that
the defined Fibonacci sequence is stationary, we have proved some results
regarding the arithmetic of integer quaternions defined on some special
division quaternion algebras and we have defined and gave properties of some
special quaternions by using Fibonacci sequences. \textbf{\ }%
\begin{equation*}
\end{equation*}

\textbf{References}%
\begin{equation*}
\end{equation*}%
\newline
\newline
[Co, Sm; 03] J.H. Conway, D.A. Smith, \textit{On Quaternions and Octonions},
A.K. Peters, Natick, Massachusetts, 2003.\newline
[Da, Sa, Va; 03] G. Davidoff, P. Sarnak, A. Valette, \textit{Elementary
Number Theory, Group Theory, and Ramanujan Graphs}, Cambridge University
Press, 2003.\newline
\textbf{[}Fib.\textbf{]}
http://www.maths.surrey.ac.uk/hosted-sites/R.Knott/Fibonacci/fib.html\newline
[Fl; 20] C. Flaut, Some applications of MV algebras,
https://arxiv.org/pdf/2001.03971.pdf\newline
\textbf{[}Fl, Sa; 15\textbf{]} C. Flaut, D. Savin, \textit{Quaternion
Algebras and Generalized Fibonacci-Lucas Quaternions}, Adv. Appl. Clifford
Algebras, 25(4)(2015), 853-862.\newline
[FSZ; 20] C. Flaut, D. Savin, G. Zaharia, \textit{Properties and
applications of some special integer number sequences}, accepted in
Mathematical Methods in the Applied Sciences, DOI: 10.1002/mma.6257\newline
(https://onlinelibrary.wiley.com/doi/abs/10.1002/mma.6257)\newline
\textbf{[}Fl, Sh; 13\textbf{]} C. Flaut, V. Shpakivskyi, \textit{On
Generalized Fibonacci Quaternions and Fibonacci-Narayana Quaternions}, Adv.
Appl. Clifford Algebras, 23(3)(2013), 673-688.\newline
[GF; 10] F. Ghaboussi, J. Freudenberger, \textit{Codes over Gaussian integer
rings}, 18th Telecommunications forum TELFOR 2010, 662-665.\newline
\textbf{[}Gi, Sz; 06\textbf{]} P. Gille, T., Szamuely, Central Simple
Algebras and Galois Cohomology, Cambridge University Press, 2006.\newline
[Gu; 13] M. G\"{u}zeltepe, \textit{Codes over Hurwitz integers}, Discrete
Math., \textbf{313(5)(2013)}, 704-714.\newline
\textbf{[}Ho; 63\textbf{]} A. F. Horadam, \textit{Complex Fibonacci Numbers
and Fibonacci Quaternions}, Amer. Math. Monthly, 70(1963), 289--291.\newline
[Hu; 94] K. Huber, \textit{Codes over Gaussian integers}, IEEE Trans.
Inform. Theory, 40(1994)\textbf{, } 207--216.\newline
\textbf{[}Ko\textbf{]} D. R. Kohel, \textit{Quaternion algebras}, available
online at\newline
http://www.i2m.univ-amu.fr/perso/david.kohel/alg/doc/AlgQuat.pdf.\newline
\textbf{[}Ma; 07\textbf{]} K. Martin, Classification of quaternion algebras,%
\newline
http://www2.math.ou.edu/~kmartin/quaint/ch5.pdf, 2017.\newline
\textbf{[}Mag\textbf{]} http://magma.maths.usyd.edu.au/magma/handbook/ 
\newline
\textbf{[}Ra; 1917\textbf{]} S. Ramanujan, \textit{On the expression of a
number in the form} $ax^{2}+by^{2}+cz^{2}+du^{2},$ Proc. Cambridge Philos.
Soc. 19(1917), 11-21. Collected Papers, AMS Chelsea Publ., Providence, RI,
2000, 169-178.\newline
[Ri; 95] J. Rif\`{a}, \textit{Groups of Complex Integer Used as QAM Signals}%
, IEEE Transactions on Information Theory, \textbf{41(5)(1995)}, 1512-1517.%
\newline
\textbf{[}Sa; 14\textbf{]} D. Savin, \textit{About some split central simple
algebras}, An. St. Univ. Ovidius Constanta Mat. Ser. 22(1) (2014), 263\^{a}%
\euro \textquotedblleft 272.\newline
\textbf{[}Sa; 17\textbf{]} D. Savin, \textit{About special elements in
quaternion algebras over finite fields} Adv. Appl. Clifford Algebras 27(2)
(2017), 1801\^{a}\euro \textquotedblleft -1813.\newline
\textbf{[}Sa; 19\textbf{]} D. Savin, \textit{Special numbers, special
quaternions and special symbol elements}, chapter in the book Models and
Theories in Social Systems, vol. 179, Springer 2019, ISBN-978-3-030-00083-7,
p. 417-430.\newline
\textbf{[}Si; 70\textbf{]} W. Sierpinski, \textit{250 Problems in elementary
number theory}, American Elsevier Pub. Co (1970).\newline
\textbf{[}Vi; 80\textbf{]} M. F. Vigneras, \textit{Arithmetique des algebres
de quaternions}, {Lecture Notes in Math., vol. 800, Springer, Berlin, 1980.}

\begin{equation*}
\end{equation*}

Cristina FLAUT

{\small Faculty of Mathematics and Computer Science,}

{\small Ovidius University of Constan\c{t}a, Rom\^{a}nia,}

{\small Bd. Mamaia 124, 900527,}

{\small http://www.univ-ovidius.ro/math/}

{\small e-mail: cflaut@univ-ovidius.ro; cristina\_flaut@yahoo.com}%
\begin{equation*}
\end{equation*}%
\medskip \qquad\ \qquad\ \ 

Diana SAVIN

{\small Faculty of Mathematics and Computer Science, }

{\small Ovidius University of Constan\c{t}a, Rom\^{a}nia, }

{\small Bd. Mamaia 124, 900527, }

{\small http://www.univ-ovidius.ro/math/}

{\small e-mail: \ savin.diana@univ-ovidius.ro, \ dianet72@yahoo.com}

\end{document}